\documentclass[12pt]{article} \textheight 8in\textwidth
6in \oddsidemargin 0in\evensidemargin 0in
\usepackage{graphicx}
\usepackage{amsfonts}
\usepackage{epsf}

\newtheorem{theorem}{Theorem}
\newtheorem{lemma}[theorem]{Lemma}

\newtheorem{conjecture}[theorem]{Conjecture}

\newtheorem{proposition}[theorem]{Proposition}
\newtheorem{problem}[theorem]{Problem}

\newcommand{\text}[1]{\quad\mbox{#1}\quad}
\def\beq{\begin{equation}}\def\eeq{\end{equation}}
\def\beqn{\begin{eqnarray}}\def\eeqn{\end{eqnarray}}
\def\pont{\hspace{-6pt}{\bf.\ }}

\def\qed{\ifhmode\unskip\nobreak\fi\quad\ifmmode\Box\else$\Box$\fi}

\title{Large cross-free sets in Steiner triple systems}

\author{Andr\'as Gy\'arf\'as  \thanks{Supported in part by
OTKA K104373.}\\[-0.8ex]\\[-0.8ex]
\small Alfr\'ed R\'enyi Institute of Mathematics\\[-0.8ex]
\small Hungarian Academy of Sciences\\[-0.8ex]
\small Budapest, P.O. Box 127\\[-0.8ex]
\small Budapest, Hungary, H-1364 \small
\texttt{gyarfas.andras@renyi.mta.hu}}
\begin{document}
\maketitle
\begin{abstract}
A {\em cross-free} set of size $m$ in a Steiner triple system $(V,{\cal{B}})$ is three pairwise disjoint $m$-element subsets
$X_1,X_2,X_3\subset V$ such that no $B\in {\cal{B}}$ intersects all the three $X_i$-s. We conjecture that for every admissible
$n$ there is an STS$(n)$ with a cross-free set of size $\lfloor{n-3\over 3}\rfloor$ which if true, is best possible. We prove this
conjecture for the case $n=18k+3$, constructing an STS$(18k+3)$ containing a cross-free set of size $6k$. We note that some of the
$3$-bichromatic STSs, constructed by Colbourn, Dinitz and Rosa, have cross-free sets of size close to $6k$ (but cannot have size exactly $6k$).

The constructed STS$(18k+3)$ shows that equality is possible for $n=18k+3$ in the following result: in every $3$-coloring of the blocks of any Steiner triple system STS$(n)$ there is a monochromatic connected component of size at least $\lceil{2n\over 3}\rceil+1$ (we conjecture that
equality holds for every admissible $n$).

The analogue problem can be asked for $r$-colorings as well, if $r-1 \equiv 1,3 \mbox{ (mod 6)}$ and $r-1$ is a prime power, we show that the answer is the same as in case of complete graphs:  in every $r$-coloring of the blocks of any STS$(n)$, there is a monochromatic connected component with at least ${n\over r-1}$ points, and this is sharp for infinitely many $n$.

\end{abstract}

\section{Introduction}

A hyperwalk in a hypergraph $H=(V,E)$ is a sequence $v_1,e_1,v_2,e_2,\dots ,v_{t-1},e_{t-1},v_t$ of vertices and edges such that for all $1\le i<t$ we have $v_i\in e_i$ and $v_{i+1}\in e_i$. We say that $v\sim w$, if there is a hyperwalk from $v$ to $w$. The relation $\sim$ is an equivalence relation, and the subhypergraphs induced by its classes are called the \textit{connected components} of $H$. A vertex $v$ that is not covered by any edge forms a trivial component with one vertex $v$ and no edge.

The size of the largest monochromatic component in edge colorings of complete graphs and hypergraphs is well investigated, for a present survey see \cite{GY1}. For example, in every $3$-coloring of the edges of the complete graph $K_n$ there is a monochromatic connected component of size at least $n/2$ and in every $3$-coloring of the edges of $K_n^3$, the complete $3$-uniform hypergraph, there is a monochromatic spanning component. What happens in between, when the blocks of a Steiner triple system $(V,{\cal{B}})$ are colored?  For example, in every coloring of the blocks of STS$(9)$ with 3 colors, there is a monochromatic connected component of size at least $7$ but in the $4$-coloring of its blocks defined by the four parallel classes, every component in every color has only $3$ points.  Let $f(n)$ denote the largest $m$ such that in every $3$-coloring of the blocks of any STS$(n)$ there is a monochromatic connected component on at least $m$ points. It is easy to see that $f(7)=6, f(9)=7$. Our main result is the following.

\begin{theorem} \pont \label{main} $f(6k+3)\ge 4k+3$ with equality if $k$ is divisible by 3. Moreover, $f(6k+1)\ge 4k+2$.
\end{theorem}

In fact, the inequalities of Theorem \ref{main} are probably always sharp (one can easily check the cases $k=1,2$):

\begin{conjecture}\pont \label{conj1} For every positive integer $k$, $f(6k+1)=4k+2, f(6k+3)=4k+3$.
\end{conjecture}

Three pairwise disjoint $m$-element sets of points, $X_1,X_2,X_3\subset V$, is a {\em cross-free set of size $m$} in a Steiner triple system $(V,{\cal{B}})$ if no block $B\in {\cal{B}}$ intersects each $X_i$ in exactly one point. To obtain the upper bound in Theorem \ref{main}, we need some STS with a cross-free set of size almost $n/3$. In Theorem \ref{const} we construct an STS$(6k+3)$ for the case $k\equiv 0 \mbox{ (mod 3)}$ which contains a cross-free set of size $2k$ (and this is best possible).

It is worth noting that constructions of Colbourn, Dinitz and Rosa in \cite{CDR} provides STS$(n)$-s with cross-free sets of size asymptotic to $n/3$. They construct $3$-bichromatic STSs where all points are partitioned into $X_1,X_2,X_3$ so that every block intersects {\em precisely two} of the $X_i$-s and they can also control the sizes of the $X_i$s. In particular, they provide $3$-bichromatic STS(n)s where the sizes are nearly equal to $n/3$. However, it follows easily that in $3$-bichromatic STS$(n)$s with $|X_1|\le |X_2|\le |X_3|$, ${n/3}-|X_i|$ tends to infinity with $n$. Therefore, to achieve a cross-free set of size $2k$ in an STS$(6k+3)$ the number of blocks inside the $X_i$s tends to infinity with $n$.

To see the connection of cross-free sets to $f(n)$, let  $G(n)$ be the size of the largest cross-free set present in {\bf some} STS$(n)$.

\begin{lemma}\pont \label{gn} $f(n)\le n-G(n)$.
\end{lemma}
\noindent {\bf Proof. } Suppose $|X_1|=|X_2|=|X_3|=G(n)$ for a cross-free set $X_1,X_2,X_3\subset V$ in some STS$(n)$. Then coloring any block $B$ with the smallest $i$ such that $B\subset V\setminus X_i$, we have a $3$-coloring of the blocks with one nontrivial monochromatic connected component of size $n-G(n)$ in each color. \qed

The next result implies the equality $f(6k+3)=4k+3$ for $k$ divisible by $3$ in Theorem \ref{main}.
\begin{theorem} \pont \label{const} For $n=18k+3$, $G(n)=6k$.
\end{theorem}

In fact, Theorem \ref{const} probably can be extended, it would imply Conjecture \ref{conj1}

\begin{conjecture}\pont \label{conj2} $G(6k+3)=2k, G(6k+1)=2k-1$.
\end{conjecture}

It is easy to see that Conjecture \ref{conj2} is sharp (if true). Indeed, a cross-free set of size $2k+1$ in an STS$(6k+3)$ would mean that there are at most $3{2k+1\choose 2}$ blocks and that is less than ${6k+3\choose 2}/3$. Similarly, a cross-free set of size $2k$ in an STS$(6k+1)$ would show that there are at most $3k+3{2k\choose 2}$ blocks, less than ${6k+1\choose 2}/3$.

One can define $f_r(n)$ similarly for $r$-colorings of blocks. A lower bound on it can be easily derived from known results.

\begin{proposition}\pont \label{lboundr} $f_r(n)\ge \lceil {n\over r-1} \rceil.$
\end{proposition}

\noindent {\bf Proof. } Any $r$-coloring of the blocks of an STS$(n)$ defines an $r$-coloring of the edges of $K_n$, by coloring the three pairs defined by a block with the color of the block. In this coloring there is a monochromatic, say red connected component $C$ with at least $\lceil {n\over r-1} \rceil$ vertices, proved first in \cite{GY}, a more accessible account is the survey \cite{GY1}. The blocks covering the red edges of $C$ obviously span a red connected component on $C$. \qed

The lower bound of Proposition \ref{lboundr} is trivially sharp for $r=2$ but also for certain other values of $r$, starting with $r=4,8,10,14,..$.

\begin{proposition}\pont \label{rcolor} $f_r(n)={n\over r-1}$ for infinitely many $n$ if $r-1$ is in the form $3^m,p^m,q^{2m}$ where $m\ge 1$, $p,q$ are primes,  $p\equiv 1 \mbox{ (mod 6)}, q\equiv -1 \mbox{ (mod 6)}$.
\end{proposition}

\noindent {\bf Proof. } $f_r(n)\ge {n\over r-1}$ follows from Proposition \ref{lboundr}. Suppose $r-1$ is a prime power in the form $3^m,p^m,q^{2m}$ where $m\ge 1$, $p,q$ are primes,  $p\equiv 1 \mbox{ (mod 6)}, q\equiv -1 \mbox{ (mod 6)}$.  This implies that $r-1\equiv 1 \mbox{ (mod 6)}$ or $r-1\equiv 3 \mbox{ (mod 6)}$. Then there exists an affine plane $P$ of order $r-1$ and we can define an STS$((r-1)^2)$ by substituting each line of $P$ by a copy of  an STS$(r-1)$. Then the blocks of STS$((r-1)^2)$ can be naturally colored with $r$ colors according to the $r$ parallel classes of $P$. In this coloring every component has size $r-1={(r-1)^2\over r-1}$, providing an example with equality. To get infinitely many, we can apply the well-known direct product construction (see \cite{CR}) of STS$(n_1n_2)$ from STS$(n_1)$,STS$(n_2)$.  Assume we already know that for some $t\ge 0$ the blocks of STS$(3^t(r-1)^2)$ can be $r$-colored so that each color class has $r-1$ nontrivial components (of size $3^t(r-1)$) and consider the STS$(3^{t+1}(r-1)^2)$ defined as  STS$(3^t(r-1)^2)$  $\times T$ where $T$ is a single block on three points. Then each component $C$ in each color class of STS$(3^t(r-1)^2)$  defines a component $C\times T$ in STS$(3^{t+1}(r-1)^2)$  whose blocks in $C\times T$ can be colored with the same color. This defines a natural $r$-coloring of the blocks of STS$(3^{t+1}(r-1)^2)$, preserving the property that each color class has $r-1$ nontrivial components. \qed

Our problem to determine $f(n)$ led to find $G(n)$, the size of the largest cross-free set present in {\bf some} STS$(n)$.
It seems natural and interesting to find or estimate the size $g(n)$ of the largest cross-free set present in {\bf any} STS$(n)$.  Obviously, $$G(n)\ge g(n)\ge {\alpha(n)\over 3}$$ where $\alpha(n)$ is the largest independent set present in {\em any} STS$(n)$. For the most recent result and history on $\alpha(n)$ see \cite{EV}.

\begin{problem} \pont Is $g(n)$ significantly smaller than $G(n)$?
\end{problem}

\section{Proof of Theorems \ref{main}, \ref{const}}

We prove first that $f(6k+3)\ge 4k+3, f(6k+1)\ge 4k+2$.

Suppose that the blocks of an STS $(V,{\cal{B}})$ with $|V|=n$ are 3-colored and consider the three components $C_1,C_2,C_3$ in colors $1,2,3$ containing a point $v\in V$. There are some cases according to the number of $C_i$s with points covered only by $C_i$, we call such points as ``private parts'' of  $C_i$.

\noindent {\bf Case 1.} No $C_i$ has private part. In this case the sets $C_i$ doubly cover $V\setminus \{v\}$ and $v$ is triply covered. This implies easily that $f(6k+3)\ge 4k+3$ and also $f(6k+1)\ge 4k+2$, unless if the $C_i$s intersect in one point and all the three doubly covered sets have size $2k$. However, in this case we can have  $3k$ blocks covering $v$ and any other block must cover a pair of $(C_i\cap C_j)\setminus \{v\}$. Thus altogether we have at most $3k+3{2k\choose 2}<{{6k+1\choose 2}\over 3}$ blocks in STS$(6k+1)$ and that is a contradiction.

\noindent {\bf Case 2.} Only $C_1$ has a private part. Now there is no point $w\in V$ that belongs to $(C_2\cap C_3)\setminus C_1$, otherwise no block can cover $wx$ where $x$ is from the private part of $C_1$. Thus in this case $C_1$ covers $V$.

\noindent {\bf Case 3.} Two $C_i$s, say $C_1,C_2$ have private parts. Now $(C_1\cap C_3)\setminus C_2$ and $(C_1\cap C_2)\setminus C_3$ are both empty and any pair of points $x,y$ from the private parts of $C_1,C_2$, respectively, must be in a block colored with color $3$. Thus the union of the private parts of $C_1,C_2$ is part of a component $C$ of color $3$. We can now apply the argument of Case 1 to the components $C,C_1,C_2$.

\noindent {\bf Case 4.} All $C_i$s have private parts. Now sets covered by precisely two of $C_1,C_2,C_3$ must be empty and the private parts $X_i\subset C_i$ together with $X_4=C_1\cap C_2\cap C_3$, partition $V$.  Pairs of points $x\in X_1,y\in X_2$ must be in a block of color $3$, pairs of points $x\in X_1,y\in X_3$ must be in a block of color $2$, pairs of points $x\in X_2,y\in X_3$ must be in a block of color $1$, thus the union of any two $X_i$s is covered by (in fact equal to) a monochromatic component.  Observe that every block of our $(V,{\cal{B}})$  must contain a pair from some of the $X_i$s, thus
\begin{equation}\label{a}s=\sum_{i=1}^4 {|X_i|\choose 2}\ge {{n\choose 2}\over 3}.
\end{equation}
First let $n=6k+3$, assume w.l.o.g that $$|X_1|\le |X_2|\le |X_3|\le |X_4|.$$
If $|X_1|\ge 3k+1+t$ for some positive integer $t$ then let $X_j$ be the largest among $X_2,X_3,X_4$. Then
$$|X_1|+|X_j|\ge 3k+1+t+{3k-t+2\over 3}\ge 4k+3$$
proving what we need.  However, if $|X_1|\le 3k+1$ then the maximum of $s$ (under the condition that each component has size at most $2k+2$) is obtained when $|X_1|=3k+1, |X_2|=|X_3|=k+1, |X_4|=k$. But this contradicts to (\ref{a}). Similar argument works if $n=6k+1$, then
$$|X_1|=3k+1, |X_2|=|X_3|=|X_4|=k$$ gives the largest $s$ and the contradiction.

This finishes the proof of the two inequalities of Theorem \ref{main}. It is left to prove that $f(6k+3)= 4k+3$ if $k$ is divisible by 3, i.e. to prove Theorem \ref{const}. In fact we need to prove only that $G(n)\ge 6k$, however $G(n)<6k+1$ follows easily: a partition of $V$ for a STS$(18k+3)$ into three sets of size $6k+1$ cannot be cross-free since then there are at most $t=3{6k+1\choose 2}$ blocks and $t$ is less than the number of blocks required in an STS$(18k+3)$.

We construct an STS$(18k+3)$ with a cross-free set of size $6k$ as follows. Let $H_k$ be the graph with $6k$ vertices and $4k$ edges, having $2k$ components, $k$ of them a $P_4$, a path on four vertices, and $k$ of them a single edge. We call the {\em middle } of $H_k$ the union of the middle edges of the $P_4$ components in $H_k$.   A {\em near factor} of a graph with $2m$ (or $2m-1$) vertices means $m-1$ pairwise disjoint edges.

\begin{lemma}\pont\label{part} Let $T$ be the graph containing $k$ vertex disjoint edges on $6k$ vertices. Then the edge set of $G_k=K_{6k}\setminus T$ can be partitioned into $2k$ factors $F_1,\dots,F_{2k}$ and $4k$ near factors $E_1,\dots,E_{4k}$ in such a way that the  pairs uncovered by the near factors form a graph isomorphic to $H_k$ and in the isomorphism the middle of $H_k$ corresponds to the pairs of $T$.
\end{lemma}

Based on the lemma, we define an STS$(18k+3)$ with a cross-free set of size $6k$. Take three disjoint copies of $H_k$ on vertex sets $X_0,X_1,X_2$ and define $\cal{T}$ as a partial triple system  PTS$(18k)$ on $\cup_{i=0}^2 X_i$ as follows. Partition each $X_i$ into $k$ $P_4$ paths $a_{6j+1}^i,a_{6j+2}^i,a_{6j+3}^i,a_{6j+4}^i$ and $k$ edges $a_{6j+5}^i,a_{6j+6}^i$ for $j=0,\dots k-1$.
This way each $X_i$ spans a copy of $H_k$.

Now Lemma \ref{part} can be applied with vertex set $X_0$ to obtain $2k$ factors and $4k$ near factors with the required properties (with respect to the copy of $H_k\subset X_0$). We can extend these factors and near-factors to blocks of $\cal{T}$, using vertices of $X_1$ as follows. Let $a_{6j+4}^1$ define blocks with the pairs of the near factor $E_{4j+1}$ with uncovered pair $a_{6j+2}^0,a_{6j+3}^0$, $j=0,\dots, k-1$. Then $a_{6j+1}^1$ defines blocks with the pairs of the near factor $E_{4j+2}$ with uncovered pair $a_{6j+5}^0,a_{6j+6}^0$, $j=0,\dots, k-1$; similarly $a_{6j+6}^1$ defines blocks with the pairs of the near factor $E_{4j+3}$ with uncovered pair $a_{6j+1}^0,a_{6j+2}^0$, $j=0,\dots, k-1$; and $a_{6j+5}^1$ defines blocks with the pairs of the near factor $E_{4j+4}$ with uncovered pair $a_{6j+3}^0,a_{6j+4}^0$, $j=0,\dots, k-1$. Finally, $a_{6j+2}^1,a_{6j+3}^1$ define blocks with the pairs of the factors $F_{2j+1},F_{2j+2}$, $j=0,\dots,k-1$.

The construction of the previous paragraph can be repeated cyclically, defining blocks with one vertex in $X_2$ and two in $X_1$, and a third time defining blocks with one vertex in $X_0$ and two in $X_2$. By Lemma \ref{part}, the partial STS $\cal{T}$ defined this way covers all pairs of $X_0\cup X_1\cup X_2$ except a $3$-regular graph $U$ of with the following edges: $a_{6j+2}^i,a_{6j+3}^i$ for $i=0,1,2$ and $j=0,\dots, k-1$ (formed by the middle of the three copies of $H_k$) and the $3\times 8k$ edges between the pairs $X_i,X_j$ that belong to the uncovered pairs of the $3\times 4k$ near factors. It can be easily seen that the graph $U$ can be factored into three $1$-factors. In fact, these factors are
$$a_{6j+2}^i,a_{6j+3}^{i},a_{6j+1}^i,a_{6j+6}^{i-1},a_{6j+5}^i,a_{6j+4}^{i-1},$$
$$a_{6j+4}^i,a_{6j+2}^{i-1},a_{6j+5}^i,a_{6j+3}^{i-1},a_{6j+6}^i,a_{6j+1}^{i-1},$$
$$a_{6j+1}^i,a_{6j+5}^{i-1},a_{6j+4}^i,a_{6j+3}^{i-1},a_{6j+6}^i,a_{6j+2}^{i-1},$$
where $i=0,1,2$ and $j=0,\dots k-1$ with arithmetic on $i,j$-s are modulo $3,6$, respectively.

Finally, $\cal{T}$ is extended to an STS$(18k+3)$ by extending each factor of $U$ to a block with one of three new points $A,B,C$ which also forms the last block. This finishes the proof of Theorem \ref{const} and with it Theorem \ref{main}. \qed

\bigskip

\noindent {\bf Proof of Lemma \ref{part}.} The required partition is constructed from the {\em standard factorization} of $K_{6k}$ on vertex set $\{1,2,\dots,6k-1\}\cup \infty$ where (for $i=1,2,\dots,6k-1$) factor $F_i$ contains $(i,\infty)$ and $\{(i-j,i+j): 1\le j \le 3k-1\}$ with mod $6k-1$ arithmetic.

We shall keep $2k-1$ of the factors $F_i$ and define the near factors $E_1,\dots,E_{4k}$ by deleting  one edge from each of the other $4k$ factors so that the deleted edges form a graph isomorphic to $H_k$. The factor formed by the middle of $H_k$ is left uncovered and all other edges of $H_k$ form a new factor $F^*$ which is added as the $2k$-th factor in the partition. To define the construction, it is enough to specify the set of $4k$ pairs (all from different $F_i$) which form a graph $Z_k$ isomorphic to $H_k$. The construction is the simplest for $k\equiv 1 \pmod 2$ so we describe that first.

Suppose that $k\equiv 1 \pmod 2$ and set $W=\{(1,3),(2,4),(3,5),(5,\infty)\}$. Moreover, for $m\in \{6,12,\dots,6(k-2)\}$ let $L_m=A_m\cup B_m$ be the following set of eight pairs on twelve consecutive numbers: $$A_m=\{(m,m+2),(m+2,m+4),(m+4,m+6),(m+1,m+3)\},$$ $$B_m=\{(m+5,m+7),(m+7,m+9),(m+9,m+11),(m+8,m+10)\}.$$
It is immediate to check that $W,A_m,B_m$ are all define ($6$-vertex) graphs with a $P_4$ component and a $K_2$ component. Thus the graph $Z_k$ defined by $W$ (for $k=1$) and by $W\cup_{m=6}^{6(k-2)} L_m$ (for odd $k>1$) is isomorphic to $H_k$. Moreover, since all edges (apart from $(5,\infty)$) of $Z_k$ are in the form $(j,j+2)$ and $j\ne 4$, each edge of $Z_k$ belongs to different $F_i$.

The case $k\equiv 0 \pmod 2$ is slightly more involved, we use another type of components $C_m,D_m$ (beside $W$) to define $Z_k$.

$$C_m=\{(m,m+1),(m,m+2),(m+2,m+4),(m+3,m+5)\},$$ $$D_m=\{(m,m+2),(m+1,m+2),(m+1,m+3),(m+4,m+5)\}.$$
For $k=2$ we use $W$ followed by $C_6$ to define $Z_2$. For $k>2$ start with $W$, then ${k\over 2}$ copies of $C_m$ ($m=6,12,\dots,3k$) then ${k-2\over 2}$ copies of $D_m$ ($m=3k+6,\dots,6(k-1)$). To check here that each edge of $Z_k$ belongs to different $F_i$, note that ``jumping pairs'' $(j,j+2)$ are obviously from different $F_i$ (from $F_{j+1}$). The same is true for the ``consecutive pairs'' $(j,j+1)$. To check consecutive pairs against jumping pairs, notice that for $m=6,12,\dots,3k$ the pair $(m,m+1)$ of $C_m$ belongs to $F_{3k+m}$, a starting point of the $D$-block opposite to $C_m$ thus it is not skipped by any jumping pair. Similarly, for $m=3k+6,\dots,6(k-1)$, the pairs $(m+1,m+2)$ and $(m+4,m+5)$ in $D_m$ belong to $F_{m+2-3k}$ and $F_{m+5-3k}$, respectively, and they are not skipped in their opposite $C$-blocks.  \qed

\noindent {\bf Acknowledgement. } Thanks for the remarks of the referees, especially for calling my attention to \cite{CDR}, showing the connection of $3$-bichromatic Steiner triple systems and the ones with large cross-free sets.

\end{document}